\ifx\control\undefined

\documentclass[10pt, leqno, twoside]{article}
\usepackage{amssymb, amsmath, bm} 
\usepackage{latexsym} 
\usepackage{graphics} 
\usepackage{url} 
\usepackage{enumerate}
\usepackage[all]{xy}
\usepackage{graphicx}
\usepackage{float}
\usepackage{amsrefs}
\usepackage{color}

\usepackage{authblk}
\usepackage{blindtext}

\usepackage{enumitem}

\numberwithin{equation}{section}

\usepackage{caption}
\usepackage{subcaption}

\newtheorem{theo}{Theorem}

\newtheorem{lemm}[theo]{Lemma}

\newenvironment{proo}[1][\proofname]{\normalfont{\itshape
#1{:}}\quad\mdseries\ignorespaces}
{{$\Box$}{\vskip\belowdisplayskip}}
\newcommand{\proofname}{Proof}
\newtheorem{defi}[theo]{Definition}

\begin{document}

\fi
%%%%%%%%%%%%%%%%%%%%%%%%%%%%%%%%%%%%%%%%%%%%%%%%%%%%%%

\title{The concept of center as an equivariant map and a proof of an analogue of the center conjecture for equifacetal simplices}

%\author{L. Felipe Prieto-Mart\'inez, Raquel S\'anchez-Cauce}

\author{L. Felipe Prieto-Mart\'inez\thanks{Department of Applied Mathematics, Universidad Polit\'ecnica de Madrid (Spain), luisfelipe.prieto@upm.es}}% \, and\, Raquel S\'anchez-Cauce\thanks{Department of Artificial Intelligence, Universidad Nacional de Educaci\'on a Distancia (Spain), rsanchez@dia.uned.es}}
%\author{}
%\affil[1]{Department of Mathematics, Universidad Aut\'onoma de Madrid (Spain)}
%\affil[2]{Department of Artificial Intelligence, Universidad Nacional de Educaci\'on a Distancia (Spain)}

\date{\today}

\maketitle

\maketitle

\begin{abstract} Several authors have remarked the convenience of understanding the different notions of center appearing in Geometry (centroid of a set of points, incenter of a triangle, center of a conic and many others) as functions. The most general way to do so is to define centers as equivariant maps between $G$-spaces. In this paper, we prove that, under certain hypothesis, for any two $G$-spaces $\mathcal A,\mathcal X$, for every $V\in \mathcal A$ and for every point $P\in\mathcal X$ fixed by the symmetry group of $V$, there exists some equivariant map $\mathfrak Z:\mathcal A\to \mathcal X$ such that $\mathfrak Z(V)=P$. As a consequence of this fact, we prove an analogue (for non-neccessarily continuous centers) of the \emph{center conjecture for equifacetal simplices}, proposed by A. L. Edmonds.

\medskip


%\medskip

%\noindent \textbf{Mathematics Subject Classification (2020):} Primary 51M04 ,  Secondary 51M15
\end{abstract}

\section{Introduction}\label{sec1}

%\emph{Long ago, someone drew a triangle and three segments across it. Each segment started at a vertex and stopped at the midpoint of the opposite side. The segments met in a point. The person was impressed and repeated the experiment on a different shape of triangle. Again the segments met in a point. The person drew yet a third triangle, very carefully, with the same result. He told his friends. To their surprise and delight, the coincidence worked for them, too. Word spread, and the magic of the three segments was regarded as the work of a higher power. }

%\cite{W}

The term \emph{center} is used, with different meanings, in many contexts and it is of great interest in Applied Mathematics (see, for instance \cite{J, R}).

In a set of papers published in the decade of 1990, C. Kimberling explained the convenience of defining triangle centers (including centroid, incenter, orthocenter, circumcenter and many others) as functions from a set of triangles $\mathcal A$ to the plane $\mathcal X$ that commute with respect to similarities (see, for example, \cite{K1,K2}). This approach allowed him, among other things, to make a list of all known triangle centers in \cite{K}. Many authors have followed an approximation, similar to the one by Kimberling, in other contexts. For instance, we can find a definition of $n$-gons center in \cite{FP, P, PS}  and a similar one of center of a set of $n$-points in  \cite{E, P}. So, in Section \ref{section.definition} of the present paper, we push the concept of center to its maximum generality in order to provide a general formal definition in terms of equivariant maps between $G$-spaces.

In the theory of classical triangle centers, we can find the following result:

\begin{theo} \label{theo.exam} Let $\mathcal A$ define the set of all triangles in the plane. For every triangle $V\in\mathcal A$, we have that:

\begin{itemize}
\item[(a)] The incenter and the orthocenter of $V$ coincide if and only if $V$ is equilateral (a proof may be found in \cite{I}).
%\item[(b)] The incenter, the centroid and the orthocenter
%are collinear if and only if the triangle is isosceles (see \cite{F}).

\item[(b)] Let $S$ be the group of congruences that fix $V$. Then for every point $P\in\mathbb R^2$ fixed by $S$, there exists a triangle center $\mathfrak Z:\mathcal A\to\mathbb R^2$ such that $\mathfrak Z(V)=P$ (see \cite{FP} for a proof of an equivalent statement).
\end{itemize}

\end{theo}

\noindent There have been several attempts in the bibliography to generalize this theorem (or part of it). For example,

\begin{itemize}

\item in \cite{AHK} we can find similar results for quadrilaterals, 

\item in \cite{P} there is a generalization of this result for $\mathcal A$ being the family of sets of $n$-points and $\mathcal A$ being the set of $n$-gons and

\item  in \cite{E} A. L. Edmonds establishes (and proves for dimension less or equal than 6) what he calls \emph{the center conjecture for equifacetal simplices}, which is a generalization for simplices of the Statement (b) in Theorem \ref{theo.exam}.

\end{itemize}

\noindent The main theorem in the present paper (Theorem \ref{theo.sym}) appears in Section \ref{section.sym} and it is a generalization of Statement (b) in Theorem \ref{theo.exam} for equivariant maps. It states that that if two $G$-spaces $\mathcal A,\mathcal X$ satisfy certain compatibility conditions, for every $V\in \mathcal A$ and for every point $P\in\mathcal X$ fixed by the symmetry group of $V$, there exists some center $\mathfrak Z:\mathcal A\to \mathcal X$ such that $\mathfrak Z(V)=P$.

As a consequence of Theorem \ref{theo.sym} we prove an analogue of the \emph{center conjecture for equifacetal simplices} in Section \ref{section.equifacetal} (Theorem \ref{theo.equifacetal}), for centers which are non-neccessarily continuous.

%%%%%%%%%%%%%%%%%%%%%%%%%%%%%%%%%%%%%%%%%%%%%%%%%%%%%%%%%%%%%%%%%%%%%%%%%%%%
\section{Definition of center and some examples}\label{section.definition}

 Let $(G,\cdot)$ be a group whose identity element is denoted by $e$. A \textbf{(left) $G$-action} in the set $\mathcal X$ is a function $\alpha:G\times \mathcal X\to \mathcal X$, for which we will use the notation $g\cdot x$ or simply $gx$ instead of $\alpha(g,x)$, satisfying that 
$$\forall x\in \mathcal X,\; e\cdot x=x\quad\text{ and }\quad\forall g,h\in G,\;\forall x\in \mathcal X,\; g\cdot(h\cdot x)=(g\cdot h)\cdot x$$

\noindent A \textbf{$G$-set} is a set $\mathcal X$ endowed with a $G$-action. Let $\mathcal A,\mathcal X$ be two $G$-sets. An \textbf{equivariant map} is a function $\mathfrak Z:\mathcal A\to \mathcal X$ that commutes with respect to the $G$-action, that is, $\forall g\in G$, $\forall V\in\mathcal A$, $\mathfrak Z(g\cdot V)=g\cdot \mathfrak Z(V)$.

In view of the literature, there seems to be a difference between the geometric concept of center and the general concept of equivariant map. Let $\mathcal A,\mathcal X$ be two $G$-sets. If $\mathcal A$ is family of subsets of $\mathcal X$, we say that the action of $G$ on $\mathcal A$ is \textbf{point-wise} if it satisfies
$$g\cdot S=\{g\cdot x:x\in S\}. $$

\begin{defi} A \textbf{center} is an equivariant map $\mathfrak Z:\mathcal A\to \mathcal X$ where $\mathcal A$ is a family of subsets of $\mathcal X$ and the action of $G$ in $\mathcal A$ is point-wise.

\end{defi}

\medskip

\noindent \textbf{Remark} \emph{We may commit two abuses of notation when working with centers. First, we admit $\mathcal A$ to be a family of subsets of $\mathcal X$ with some extra structure, such as sequences, polygons or graphs. In this paper, this abuse will only appear in Example (IV) in this section. Second, for some $V\in\mathcal A$ and for some center $\mathfrak Z:A\to X$, we use the term \emph{center} of $V$ to refer to $\mathfrak Z(V)$.}

\medskip

%Let $\mathcal F$ be a family of subset of $X$ which is \textbf{closed with respect to the $G$-action}, that is, if $A\in\mathcal F$ then,  $\forall g\in G$, $g(A)\in\mathcal F$

In the following we may find some examples of centers appearing in Geometry:

\begin{itemize}

\item[(I)] The most important examples  (for us)  of equivariant maps are triangle centers (for the sake of uniformity we are not going to present the same definition as C. Kimberling). In this case, the set $\mathcal X$ is $\mathbb R^2$ and $\mathcal A$ is a family of sets with three elements in $\mathcal X$ (representing the vertices). Elements in $\mathcal A$ are called \textbf{triangles}. $G$ is the group of similarities in $\mathcal X$ with the usual action and the action of each $g\in G$ over $\mathcal A$ is
$$g\{A,B,C\}=\{gA,gB,gC\}.$$

\noindent  The simplest example of triangle center is the centroid $\mathfrak C:\mathcal A\to\mathbb R^2$ given by
$$\mathfrak C(\{A,B,C\})=\frac{1}{3}A+\frac{1}{3}B+\frac{1}{3}C. $$

\noindent Note that for some centers the domain may fail to be the whole family of sets of three points in $\mathbb R^2$:  the circumcenter (with its usual definition as the function that maps each triangle to the point where the perpendicular bisectors of the sides meet) is not defined if the three vertices are collinear.

\item[(II)] The centroid $\mathfrak C:\mathcal A\to \mathcal X$ can be defined for many other choices of $\mathcal X, \mathcal A$. For instance, it can be done for $\mathcal X=\mathbb R^n$ and $\mathcal A$ being the family of all finite sets in $\mathcal X$ and in this case is given by
\begin{equation}\label{eq.centroid} \mathfrak C(\{P_0,\ldots, P_n\})=\frac{1}{n+1}P_0+\ldots+\frac{1}{n+1}P_n. \end{equation}

\item[(III)] Let $\mathcal X$ be the real projective plane, $\mathcal A$ be the set of conic sections in the projective plane (not including the degenerate cases) and $G$ be the set of projective transformations of $\mathcal X$ (with the natural action over $\mathcal X$). The function $\mathfrak Z:\mathcal A\to \mathcal X$ that maps each conic to the so called \emph{center of the conic} (pole of the line at infinity) is also a center in the sense described above.

% \textcolor{red}{COMPROBAR Y CITAR}

\item[(IV)]  Let $\mathcal X=\mathbb R^2$. For $n\geq 4$, a $n$-gon is not simply a set of $n$ points, but a set $V$ consisting of $n$ points, called \emph{vertices}, together with and \emph{adjacency relation}: each vertex is adjacent to exactly two vertices and, for every $P_1,P_k\in V$, for $2\leq k\leq n$, there are exactly two chains $P_1,P_2,\ldots,P_k$ such that all the elements are different and two consecutive vertex are adjacent. Some examples of polygon centers may be found in \cite{FP, P, PS}. In this case, that corresponds to the first abuse of notation explained in the remark above, two different polygons may have the same vertices and a given center may assign different points to them. 

\end{itemize}

%\textcolor{red}{OPERACION N-ARIA PERMITE UNA DEFINICION DE CENTRO N-ARIO}

%%%%%%%%%%%%%%%%%%%%%%%%%%%%%%%%%%%%%%%%%%%%%%%%%%%%%%%%%%%%%%%%%%%%%%%
\section{A result for equivariant maps} \label{section.sym}

Let $\mathcal X$ be a $G$-set.  We say that $P\in \mathcal X$  is \textbf{fixed} by $g\in G$ (resp. by a subgroup $H<G$) if $g\cdot P=P$ (resp. if $\forall g\in H$, $g\cdot P=P$).  The \textbf{orbit} of a given element $P\in \mathcal X$ is the set $G\cdot x=\{g\cdot P:g\in G\}$. The action of $G$ in $\mathcal X$ is \textbf{transitive} if for every two points $P,Q\in\mathcal X$, there is some $g\in G$ such that $g\cdot P=Q$, that is, all the elements in $\mathcal X$ lie in the same orbit.

\medskip

\noindent \textbf{Notation} \emph{We denote by $\mathcal X^g$ (resp. $\mathcal X^H$) to the set of points in $\mathcal X$ fixed by $g$ (resp. by $H$). For any $P\in X$, we denote by $G_{P}$ to the \textbf{stabilizer} or \textbf{group of symmetries} of $P$, that is, the subgroup of $G$ consisting in those elements $g\in G$ such that $g\cdot P=P$. }

\medskip

%We say that a set $S\subset \mathcal X$ is \textbf{fixed} by $g\in G$ (resp. by $H<G$) if  for every $P\in S$, $P$ is fixed by $g$ (resp. for every $P\in\mathcal X$, $P$ is a point fixed by $H$). 

%We say that the action of $G$ is \textbf{transitive} if, for every $P,Q\in\mathcal X$, there is some $g\in G$ such that $Q=g\cdot P$.

Let us begin with the following lemma of immediate proof:

\begin{lemm} \label{lemma} Let $\mathcal A,\mathcal X$ be two $G$-spaces and $g\in G$. If there is some equivariant map $\mathfrak Z:\mathcal A\to \mathcal X$ and $V\in\mathcal A^g$, then $\mathfrak Z(V)\in\mathcal X^g$. In particular:
\begin{itemize}

\item[($\star$)] For every $g\in G$, if $\mathcal A^g$ is non-empty, then so is $\mathcal X^g$.

\end{itemize}

\end{lemm}

\begin{proo} Let $V\in\mathcal A^g$. If $\mathfrak Z$ is equivariant, we have that $\mathfrak Z(V)=\mathfrak Z(g\cdot V)=g\cdot \mathfrak Z(V)$. The statement ($\star$) is a consequence of this fact.

\hspace{12cm}\end{proo}

This simple observation leads us to the following, which is a generalization of Statement (b) in Theorem \ref{theo.exam} and of part of the results in \cite{P}:

\begin{theo} \label{theo.sym} Let $\mathcal A, \mathcal X$ be two $G$-spaces satisfying condition ($\star$). Let $V\in \mathcal A$. For every $P\in \mathcal X$ fixed by the group of symmetries of $V$, there exists some equivariant map $\mathfrak Z:A\to X$ such that $\mathfrak Z(V)=P$.

\end{theo}

%Suppose that $g$ is in the stabilizer $V$. According to Lemma \ref{lemma} then then $P$ is neccessary in the set of points fixed by $g$.

\begin{proo} Let $\mathcal A'$ be the orbit of $V$ in $\mathcal A$. Let us begin defining $\mathfrak Z$ for the elements in $\mathcal A'$. Let us impose that $\mathfrak Z(V)=P$.  For any $W\in \mathcal A'$ ($W$ may equal $V$), there exists at least one $g\in G$ such that  $g\cdot  V= W$, so let us define
$$\mathfrak Z(W)=g\cdot \mathfrak Z( V)=g\cdot P.$$

\noindent We need to check that this function is well defined. Suppose that there is another element $h\in G$ such that $h\cdot  V= W$. Then, we need that
$$\mathfrak Z( W)=g\cdot \mathfrak Z( V)=h\cdot \mathfrak Z(V)$$

\noindent Note that
$$g\cdot V=h\cdot V\Leftrightarrow V=g^{-1}\cdot h\cdot V, \qquad g\cdot \mathfrak Z( V)=h\cdot \mathfrak Z(V)\Leftrightarrow P=g^{-1}\cdot h\cdot P$$

\noindent and according to the statement of the theorem, {$V=g^{-1}\cdot h\cdot V$ implies $ P=g^{-1}\cdot h\cdot P$ and so $\mathfrak Z$ is well defined.}

If the action of $G$ over $\mathcal A$ were {transitive} (that is, $\mathcal A=\mathcal A'$) then we would have finished. In other case, to extend $\mathfrak Z$ to $\mathcal A$, we need to use the \emph{Axiom of Choice} as follows:

\begin{itemize}

\item Let $\{G_i\}_{i\in I}$ denote the set of subgroups of $G$ with a non-empty set of fixed points in $\mathcal A$. Consider also the set $\{\mathcal X^{G_i}\}_{i\in I}$. Let $\Phi_1$ be a {choice function} that selects an element $P_i$ from each $\mathcal X^{G_i}$, for $i\in I$ (here is where hypothesis ($\star$) is used).

\item $\mathcal A\setminus \mathcal A'$ can be split into a disjoint union of orbits $\mathcal A\setminus \mathcal A'=\bigcup_{j\in J}\mathcal A_j$. Let $\Phi_2$ be a choice function that selects one element $V_j$ in each orbit $\mathcal A_j$. 

\item If $V$ belongs to the orbit $\mathcal A_j$, then $V=g\cdot V_j$ for some $g\in G$. If $G_i$ is the stabilizer of $V_j$, for some $i\in I$, we can define
$$\mathfrak Z(V)=\mathfrak Z(g\cdot V_j)=g\cdot\mathfrak Z(V_j)=g\cdot P_i. $$

The same argument used before shows that this function is well defined and it is equivariant.

\end{itemize}

%\textcolor{red}{FALTAFALTAFALTA ESTE ULTIMO PARRAFO ES MEJORABLE. MEJOR PODRIAMOS COGER PARA TODOS LOS ELEMENTOS CON UN DETERMINADO ESTABILIZADOR SIEMPRE LA MISMA IMAGEN (ESO REQUIERE EL AXIOMA DE ELECCION). DESPUES PARA CADA OTRO COSO RELACIONADO CON ESTOS, LO DEFINIMOS COMO ARRIBA}

\hspace{12cm} \end{proo}

%%%%%%%%%%%%%%%%%%%%%%%%%%%%%%%%%%%%%%%%%%%%%%%%%%%%%%%%%%%%%%%%%%%%%
\section{Proof of the analogue of the center conjecture for equifacetal simplices} \label{section.equifacetal}

\noindent \textbf{Remark} \emph{Let $\mathcal A$ be the set of subsets of $\mathbb R^n$ consisting in $n+1$ points and }
$$\mathcal B=\{(P_0,\ldots,P_n)\in(\mathbb R^n)^2: P_0,\ldots,P_n\text{ are affinely independent}\} $$

\noindent \emph{In the original paper \cite{E}, centers are considered as continuous functions $\mathfrak Z:\mathcal B\to\mathbb R^n$. Here, centers will be functions $\mathfrak Z:\mathcal A\to\mathbb R^n$ (non-neccessarily continuous). From now on, the elements in $\mathcal A$ will be called }simplices \emph{and it is not required that the elements in $\mathcal A$ are affinely independent.}

%\emph{Unlike what appears in the original paper \cite{E}, we will use the term} simplex \emph{to refer to a set of $n+1$ points in $\mathbb R^n$, regardless of whether it is affinely independent or not. Moreover, let $\mathcal B$ denote the subset  the difference between the statement proved here and the original} center conjecture for equifacetal simplices \emph{proposed in \cite{E} is that here the center obtained is not a continuous function $\mathfrak Z:\mathcal B$}

\medskip

Let $E(n)$ denote the Euclidean group in $\mathbb R^n$.
  A \textbf{facet} of the simplex $V$ is a subset of $V$ consisting  of $n$ points. Recall that a simplex is called \textbf{equifacetal} if all the \emph{facets} are congruent to one another, that is, the natural action of $E(n)$ in the set of facets of the simplex is transitive.

Finally we can prove the \emph{center conjecture for equifacetal simplices} as a consequence of Theorem \ref{theo.sym}.

\begin{theo}[proof of the analogue of the center conjecture for equifacetal simplices] \label{theo.equifacetal} Let $\mathcal A$ be the family of subsets of $\mathcal X$ with $n+1$ elements ({simplices}). Consider the natural action of the Euclidean group $E(n)$ in both of them. Then all the centers of some $V\in \mathcal A$ coincide (that is, there exists some $P\in\mathcal X$ such that, for every center $\mathfrak Z:\mathcal A\to \mathbb R^n$, $\mathfrak Z(V)=P$) if and only if $V$ is equifacetal.

\end{theo}

\begin{proo} Proposition 3.4 in \cite{E} guarantees that for any equifacetal simplex, all the centers coincide.  Consider the centroid $\mathfrak C:\mathcal A\to\mathcal X$ as defined in Example (II). We only need to prove that, for every $V=\{V_0,\ldots, V_n\}\in \mathcal A$ which is not equifacetal, there is a center $\mathfrak Z:\mathcal A\to\mathcal X$ such that $\mathfrak Z(V)\neq\mathfrak C(V) $.

If $V$ is not affinely independent, let us consider $\mathfrak Z:\mathcal A\to\mathcal X$ given by
$$\mathfrak Z(\{V_0,\ldots, V_n\})=\begin{cases}\mathfrak C(V)& \text{if }V_0=\ldots=V_1\\ \frac{h_0}{h_0+\ldots+h_n}V_0+\ldots+\frac{h_n}{h_0+\ldots+h_n}V_n &\text{in other case}\end{cases}$$

\noindent where $h_i$ denotes the distance between $V_i$ and the centroid of $V\setminus\{V_i\}$. Note that $\mathfrak Z(V)=\mathfrak C(V)$  if and only if all the elements in $V$ coincide (and so $V$ is equifacetal) or
$$\frac{h_0}{h_0+\ldots+h_n}=\ldots=\frac{h_n}{h_0+\ldots+h_n}=\frac{1}{n+1}$$

\noindent But this second case is not possible since at least one of the values $h_i$ (but not all of them) equals 0.

If $V$ is affinely independent, then the fixed point set $\mathcal X^{(E(n))_V}$ of the symmetry group $(E(n))_V$ of $V$ is a single point if and
only if the group $E(n)$ acts transitively on $V$ (Lemma 3.3 in \cite{E}). In turn, $(E(n))_V$ acts transitively in $V$ if and only if it does on the set of facets (Theorem 3.3 in \cite{E1}). So if $V$ is not equifacetal, then the group of symmetry $(E(n))_V$ has more than one fixed point. Let $Q$ denote a fixed point different from $\mathfrak C(V)$. According to Theorem \ref{theo.sym} there exists some center $\mathfrak Z:\mathcal A\to\mathcal X$ such that $\mathfrak Z(V)=Q\neq \mathfrak C(V)$.

\hspace{12cm} \end{proo}

%%%%%%%%%%%%%%%%%%%%%%%%%%%%%%%%%%%%%%%%%%%%%%%%%%%%%%%%%%%%%%%%%%%%%
\section{Final comments}

It remains to prove the original \emph{center conjecture for equifacetal simplices} (concerning continuous centers). On the other hand, the \emph{strong conjecture for equifacetal simplices}, asking for a center $\mathfrak Z:\mathcal A\to\mathcal X$ satisfying the property
$$\mathfrak Z(V)=\mathfrak C(V)\Leftrightarrow V\text{ is equifacetal},$$

\noindent remains open for both, continuous and non-neccessarily continuous centers.

\end{document}